\documentclass[12pt]{article} 

\usepackage{amsmath,amsxtra,amssymb,latexsym, amscd}
\usepackage[mathscr]{eucal}

\begin{document}                        
                        \def\DATE{9/13/'03}

\title{ \huge\bf Omega-limit sets  and bounded solutions
} 

\def\1{\rule{0cm}{0cm}} \def\qd{\rule{3mm}{3mm}} \def\BB{$\bullet$}
\renewcommand{\arraystretch}{1.25}
\renewcommand{\theequation}{\thesectn.\arabic{equation}}
\def\sce{\setcounter{equation}{0}}  \newcounter{sectn} 
\newcounter{sbsect}
\def\sect#1{\addtocounter{section}{1}\sce\setcounter{sbsect}{0}%
        \renewcommand{\thesectn}{\thesection}\1\smallskip\\
        {\1\hspace{-2em}\large\bf\thesectn.\qquad #1\smallskip\par}}
\def\subsect#1{\addtocounter{sbsect}{1}\sce%
        
\renewcommand{\thesectn}{\thesection:\Alph{sbsect}}\1\smallskip\\
        {\bf\1\hspace{-1.5em}\thesectn.\qquad #1\smallskip\par}}
\newtheorem{Theorem}{THEOREM} \newtheorem{Lemma}[Theorem]{LEMMA}
\newtheorem{Corollary}[Theorem]{COROLLARY}
\def\thm#1#2{\be{Theorem}{\lb{#1} #2}} \def\LEM#1#2{\BE{Lemma}{\LB{#1} 
#2}}
\def\COR#1#2{\BE{Corollary}{\LB{#1} #2}}
\def\proof{\bigskip\noindent {\sc Proof:}\qquad}
\def\REM{\1\smallskip\par\noindent{\bf REMARK:}\qquad }
\def\qed{\hfill$\quad$\qd\medskip\\} \def\ds{\displaystyle}
\def\LB#1{\label{#1}} \def\BE#1#2{\begin{#1} #2 \end{#1}}
\def\EQ#1#2{\BE{equation}{\LB{#1} #2}} \def\ARR#1#2{\BE{array}{{#1} 
#2}}
\def\DES#1{\BE{description}{#1}} \def\QT#1{\BE{quote}{#1}}
\def\ENUM#1{\BE{enumerate}{#1}} \def\ITM#1{\BE{itemize}{#1}}
 \def\COM#1{\par\noindent{\bf COMMENT:\quad\sl #1}\par\noindent}
\def\mapsfrom{\hbox{$\;{\leftarrow}\kern-.15em{\mapstochar}\:\:$}}
\def\vv{\kern.344em{\rule[.18ex]{.075em}{1.32ex}}\kern-.344em}
\def\RE{\mbox{\rm I\kern-.21em R}} \def\CX{\mbox{\rm \vv C}}
\def\imp{\Rightarrow} \def\emb{\hookrightarrow} 
\def\wk{\rightharpoonup}
\def\rd{\dot{\1}} \def\d{\cdot} \def\+{\oplus} \def\x{\times}
\def\<{\langle} \def\>{\rangle} \def\o{\circ} \def\at#1{\Bigr|_{#1}}
\def\cd{\partial} \def\grad{\nabla} \def\L{\left} \def\R{\right}
\def\A{\mathbf  A} \def\by{\mathbf{y}} \def\bS{\mathbf{S}}
\def\I{{\cal I}} \def\A{{\mathbf A}} \def\D{{\mathbb D}}\def\bc{\mathbf{c}}
\def\X{\mathbf{X}} \def\bo{\mathbf{O}}
 \def\B{\mathbf{B}}\def\bone{\mathbf{1}}\def\bzero{\mathbf{0}}
\def\H{{\mathcal H}} \def\U{{\cal U}} \def\bc{\mathbf{c}}
\def\F{\mathbf{F}}
\def\eq{equation} \def\de{differential \eq} \def\pde{partial \de}
\def\sol{solution} \def\pb{problem} \def\bdy{boundary} 
\def\fn{function}
\def\dde{delay \de} \def\ev{eigenvalue}
\def\R{\mathbb R}
\def\Q{\mathbb Q}
\def\C{\mathbb C}
\def\N{\mathbb N}
\def\Z{\mathbb Z}
\def\X{\mathbb X}

\author{
{\bf Dang Vu Giang}\\
Hanoi Institute of Mathematics\\
Vietnam Academy of Science and Technology\\
18 Hoang Quoc Viet, 10307 Hanoi, Vietnam\\
{\footnotesize          e-mail: $\<$dangvugiang@yahoo.com$\>$}\\
\1\\
}
\maketitle
\bigskip
\noindent {\bf Abstract.}   
We prove among other things that the omega-limit set of a bounded solution  of a Hamilton system 
\[\left\{ \begin{aligned}
  & \mathbf{\dot{p}}=\frac{\partial H}{\partial \mathbf{q}} \\ 
 & \mathbf{\dot{q}}=-\frac{\partial H}{\partial \mathbf{p}} \\ 
\end{aligned} \right.\]
   is containing a full-time solution
so there are  the limits of $\frac 1t\int_0^t {\mathbf p}(s)ds$  and $\frac 1t\int_0^t {\mathbf q}(s)ds$ as $t\to\infty$ for any bounded solution $(\mathbf {p,q})$ of  the Hamilton system. These limits are stationary points of the Hamilton system so if a Hamilton system has no stationary point then every solution of this system is unbounded.



\sect{INTRODUCTION }

\bigskip
\noindent In this paper
  $\left( \mathbb{X},{{\left\| {\cdot} \right\|}_{\mathbb{X}}} \right)$ denotes a complex Banach space.
Let  $A:\X\to\X$ be a bounded linear operator with compact spectrum $ \sigma(A)$ and positive spectral radius $r(A)$.  In \cite{Giang} we proved that if $\sigma(A)\cap i \R=\left\{ {{i\xi }_{1}},{{i\xi }_{2}},\cdots ,{{i\xi }_{n}} \right\}$ then every bounded full-time  solution of differential equation $\mathbf{\dot{x}}(t)=A\mathbf{x}(t)$ 
has the form ${\mathbf u}\left( t \right)=\sum\limits_{k=1}^{n}{{{e}^{i{{\xi }_{k}t}}}{{\mathbf 
v}_{k}}}$,  
where ${{\mathbf v}_{1}},{{\mathbf v}_{2}},\cdots ,{{\mathbf v}_{n}}$ are fixed  vectors of $\mathbb{X}.$  Recall that full-time solution is the solution satisfying the differential equation for all $t\in\R$. For example,  periodic solutions (if exist)  are full-time and bounded solutions.
We used Beurling spectrum \cite{Giang}  and Fourier coefficients of a bounded function (on the real line) in the proof.
More exactly, we proved that the Beurling spectrum of any bounded full-time  solution is a subset of $\left\{ {{\xi }_{1}},{{\xi }_{2}},\cdots ,{{\xi }_{n}} \right\}.$
For the delay equation $\mathbf {\dot{u}}\left( t \right)=-\mathbf{u}\left( t-\tau  \right)$we proved that every almost periodic solution is periodic, so if there exists an almost periodic solution then the delay $\tau$ must be $\pi/2$. 
Generally, the spectrum of any bounded full-time solution of the delay equation $\mathbf{\dot{x}}(t)=A\mathbf{x}(t-\tau )$ is a compact subset of the interval $[-r(A),r(A)]$. 
Now consider a bounded solution $\mathbf{x}$ of 
\[\left\{ \begin{matrix}
   \mathbf{\dot{x}}\left( t \right)=A\mathbf{x}\left( t \right)\text{ for }t>0  \\
   \mathbf{x}\left( 0 \right)\text{ given in }\mathbb{X}.  \\
\end{matrix} \right.\] 
Assume that the orbit $\ \left\{ \mathbf{x}\left( t \right):\text{ }t\ge 0 \right\}$  is relatively compact.
Then the omega-limit set $\omega$ of $\mathbf{x}$ is a compact connected subset of $\mathbb{X}$ \cite{Hale}. Moreover, $\omega$  is invariant under the group $T(t)=e^{At}$. Let $\mathbf{v}$ be a point in this omega-limit set and $\mathbf{u}(t)=T(t)\mathbf{v}$. Then $\mathbf{u}$ is a bounded full-time  solution of the differential equation $\mathbf{\dot{x}}=A\mathbf{x}$. 
On the other hand, $\Omega =\omega \cup \left\{ \mathbf{x}\left( t \right):\text{ }t\ge 0 \right\}$  is a compact subset of $\mathbb{X}$. 
Therefore, the semi-group ${{\left\{ T\left( t \right) \right\}}_{t\ge 0}}$ acts injectively on $\Omega .$ By an ergodic theorem \cite{Sinai} we have
$\underset{t\to \infty }{\mathop{\lim }}\,\frac{1}{t}\int_{0}^{t}{\mathbf{x}}(s)ds=\underset{t\to \infty }{\mathop{\lim }}\,\frac{1}{t}\int_{0}^{t}{\mathbf{u}}(s)ds$. This limit is lying in the kernel of $A$.
Specially, if $\sigma(A)\cap i\R=\emptyset$ then 0 is the only bounded full time solution. Thus, every bounded  solution
tends to 0 as $t\to\infty$. 
Now let $(\mathbf {p,q})$  be a bounded solution of the  Hamilton system 
\[\left\{ \begin{aligned}
   \mathbf{\dot{p}}& =\frac{\partial H}{\partial \mathbf{q}} \\ 
  \mathbf{\dot{q}}& =-\frac{\partial H}{\partial \mathbf{p}}. \\ 
\end{aligned} \right.\]
Then there is an injective continuous semi-flow $T\left( t \right):{{\mathbb{R}}^{2n}}\to {{\mathbb{R}}^{2n}}$ such that $\left( \mathbf{p}\left( t \right),\mathbf{q}\left( t \right) \right)=T\left( t \right)\left( \mathbf{p}\left( 0 \right),\mathbf{q}\left( 0 \right) \right)$
Then the omega-limit set $\omega$ of $(\mathbf{p,q})$ is a compact connected subset of $\mathbb{R}^{2n}$ \cite{Hale}. Moreover, $\omega$  is invariant under the group $T(t)$. 
The dynamical system $\left\langle \omega ,{{\left\{ T(t) \right\}}_{t\in \mathbb{R}}} \right\rangle $ is uniquely ergodic, since 
the only invariant (continuous) function on $\left\langle \omega ,{{\left\{ T(t) \right\}}_{t\in \mathbb{R}}} \right\rangle $ is the constant function. 
Let $\mathbf{v}$ be a point in this omega-limit 
set and $\mathbf{u}(t)=T(t)\mathbf{v}$. Then $\mathbf{u}$ is a bounded full-time  solution of the differential equation \[\left\{ \begin{aligned}
   \mathbf{\dot{p}}& =\frac{\partial H}{\partial \mathbf{q}} \\ 
  \mathbf{\dot{q}}& =-\frac{\partial H}{\partial \mathbf{p}}.\\ 
\end{aligned} \right. \] 
 By an ergodic theorem \cite{Sinai}  there are  the limits of $\frac 1t\int_0^t {\mathbf p}(s)ds$  and $\frac 1t\int_0^t {\mathbf q}(s)ds$ as $t\to\infty$ for any bounded solution $(\mathbf {p,q})$ of  the Hamilton system. These limits are stationary points of the Hamilton system. Therefore, we have

\bigskip 
\par\noindent {\bf Theorem A. } If the gradient $\nabla H$ of a smooth hamiltonian $H$ is nowhere 0 then every solution of the  Hamilton system
\[\left\{ \begin{aligned}
  & \mathbf{\dot{p}}=\frac{\partial H}{\partial \mathbf{q}} \\ 
 & \mathbf{\dot{q}}=-\frac{\partial H}{\partial \mathbf{p}} \\ 
\end{aligned} \right.\]
 is unbounded.

 \bigskip 
\par\noindent 
For example, consider the system $\ddot x=-\sin x$ with $x(0)=0$. If $\dot x(0)>2$ then $x(t)$ is unbounded. If $\dot x(0)=2$
then 
$$x(t)=2\text{arc sin }\frac{e^{2t}-1}{e^{2t}+1}$$
which is increasingly tending to $\pi$ as $t\to\infty$. If $\dot x(0)\in (0,2)$ then $x(t)$ is periodic and bounded by $\pi$ in the time and both $\frac 1t\int\limits_0^t x(s)ds$ and $\frac 1t\int\limits_0^t\dot x(s)ds$ tend to 0 as $t\to\infty$. Moreover,  the period of this solution is 
\[2\int_0^A\frac{dx}{\sqrt{2\cos x-2+\dot x(0)^2}},
\]
where $A=\text{arc}\cos \left(1-\frac{\dot x(0)^2}2\right)$ is the maximal value of $x(t)$.

\bigskip
\bigskip
\sect{ MAIN RESULTS}

 \bigskip
\noindent
Let $T\left( t \right):\mathbb{X}\to \mathbb{X}$ for $t\ge 0$ denote a  semi-group with (unbounded and close) generator $A.$ Let $\mathbf{x}\left( t \right)=T\left( t \right)\mathbf{x}\left( 0 \right)$ denote a bounded solution of the differential equation $\mathbf{\dot{x}}=A\mathbf{x}$. 
Assume that the orbit $\ \left\{ \mathbf{x}\left( t \right):\text{ }t\ge 0 \right\}$  is relatively compact.
Then the omega-limit set $\omega$ of $\mathbf{x}$ is a compact connected subset of $\mathbb{X}$ \cite{Hale}. Moreover, $\omega$  is invariant under the semi-group ${{\left\{ T(t) \right\}}_{t\ge 0}}.$ Clearly, $T\left( t \right):\omega \to \omega $ is bijective. 
It is easy to prove that the dynamical system $\left\langle \omega ,{{\left\{ T(t) \right\}}_{t\in \mathbb{R}}} \right\rangle $ is  uniquely ergodic \cite{Sinai}. 
In fact, the only invariant (continuous) function on $\left\langle \omega ,{{\left\{ T(t) \right\}}_{t\in \mathbb{R}}} \right\rangle $ is the constant function. 
Hence, there is a unique Borel probability measure $\mu $ on $\omega$ \cite{Sinai} such that
$$\underset{t\to \infty }{\mathop{\lim }}\,\frac{1}{2t}\int_{-t}^{t}{\varphi \left( \mathbf{u}(s) \right)ds}=\int\limits_{\omega }{\varphi \left( \mathbf{v} \right)d\mu \left( \mathbf{v} \right)}.$$
Here, $\varphi $ denotes a continuous function on $\omega$ and 
$\mathbf{u}(s)=T\left( s \right)\mathbf{v}$ for some $\mathbf{v}\in \omega \mathbf{.}$
Therefore, there is the limit of $\frac{1}{t}\int_{0}^{t}{\mathbf{u}(s)ds}$ as $t\to \infty .$ Similarly, the limit of $\frac{1}{t}\int_{0}^{t}{\mathbf{x}(s)ds}$ exists as $t\to \infty .$

\bigskip 
\par\noindent {\bf Theorem B. } Let $A$ denote the generator of a linear semigroup $T\left( t \right):\mathbb{X}\to \mathbb{X}$ for $t\ge 0$. Let $\mathbf{x}\left( t \right)=T\left( t \right)\mathbf{x}\left( 0 \right)$ denote a bounded solution of the differential equation $\mathbf{\dot{x}}=A\mathbf{x}$. 
Assume that the orbit $\ \left\{ \mathbf{x}\left( t \right):\text{ }t\ge 0 \right\}$  is pre-compact.  Then the limit of $\frac{1}{t}\int_{0}^{t}{\mathbf{x}(s)ds}$ exists as $t\to \infty .$ This limit is a vector in the kernel of the operator $A$.
If $\sigma(A)\cap i \R=\left\{ {{i\xi }_{1}},{{i\xi }_{2}},\cdots ,{{i\xi }_{n}} \right\}$ then every bounded full-time  solution of differential equation $\mathbf{\dot{x}}(t)=A\mathbf{x}(t)$ 
has the form ${\mathbf u}\left( t \right)=\sum\limits_{k=1}^{n}{{{e}^{i{{\xi }_{k}t}}}{{\mathbf 
v}_{k}}}$,  
where ${{\mathbf v}_{1}},{{\mathbf v}_{2}},\cdots ,{{\mathbf v}_{n}}$ are fixed  vectors of $\mathbb{X}.$
Specially, if $\sigma(A)\cap i\R\subseteq\{0\}$ then  every bounded solution of pre-compact orbit tends to a vector in the kernel of $A$ as $t\to \infty .$

\bigskip 
\par\noindent {\sl Proof: } As we have mentioned before, the dynamics on the omega limit set of $\mathbf{x}$ is uniquely ergodic. Moreover, this limit set contains a full time bounded solution. 
Let $\mathbf {u}$ denote a bounded full-time solution of $\mathbf{\dot{x}}(t)=A\mathbf{x}(t)$. 
Then $(\lambda-D)^{-1}{\mathbf u}(t)=(\lambda-A)^{-1} {\mathbf u}(t)
 $ for any $t\in\R$ and 
$\lambda\notin i\R\cup\sigma(A)$. Here $D$ denotes the differential operator with spectrum $i\R$. Therefore, for any point $\xi$ in the Beurling spectrum of $u$ we have $i\xi\in\sigma(A)$.  Hence, if $\sigma(A)\cap i \R=\left\{ {{i\xi }_{1}},{{i\xi }_{2}},\cdots ,{{i\xi }_{n}} \right\}$  then the Beurling spectrum of any bounded full-time  solution is a subset of $\left\{ {{\xi }_{1}},{{\xi }_{2}},\cdots ,{{\xi }_{n}} \right\}.$ 
Thus, ${\mathbf u}\left( t \right)=\sum\limits_{k=1}^{n}{{{e}^{i{{\xi }_{k}t}}}{{\mathbf 
v}_{k}}}$,  
where ${{\mathbf v}_{1}},{{\mathbf v}_{2}},\cdots ,{{\mathbf v}_{n}}$ are fixed  vectors of $\mathbb{X}$
\cite{Giang}, \cite{Minh}.
Now consider a bounded solution $\mathbf x$ of pre-compact orbit. Then the omega-limit set of $\mathbf x$ 
should contain a bounded full time solution ${\mathbf u}\left( t \right)=\sum\limits_{k=1}^{n}{{{e}^{i{{\xi }_{k}t}}}{{\mathbf 
v}_{k}}}$. 
Specially, if $\sigma(A)\cap i\R\subseteq\{0\}$  then the omega limit set of any bounded solution with  pre-compact orbit has only one element. This element is a vector of the kernel of $A$.  The proof is now complete.

\bigskip 
\par\noindent {\bf Remark. } The last statement in our Theorem makes a significant extension of results in \cite{Phong1},  \cite{Phong2}.  Indeed, the authors have proved the existence of the $\lim\limits_{t\to\infty}\frac{1}{t}\int_{0}^{t}{\mathbf{x}(s)ds}$ only.

{\footnotesize
\bigskip 
\par\noindent {\bf Acknowledgement. }
 Deepest appreciation is extended towards the NAFOSTED  (the National Foundation for Science and Techology Development in Vietnam) for the financial support.

}

\end{document}